\documentclass{amsart}

\usepackage{amscd,amssymb,amsthm,mathrsfs}
\usepackage{graphicx}
\usepackage{hyperref}
\DeclareMathOperator{\RE}{Re}
 
 \newcommand{\set}[1]{\left\{#1\right\}}
\newtheorem{thm}{Theorem}[section]
\newtheorem{cor}[thm]{Corollary}
\newtheorem{exm}[thm]{Example}
\newtheorem{lem}[thm]{Lemma}

\theoremstyle{definition}
\newtheorem{defn}[thm]{Definition}
\theoremstyle{remark}
\newtheorem{rem}[thm]{Remark}

\title[On the Fekete-Szeg\"{o} inequality]{On the Fekete-Szeg\"{o} inequality for Certain class of analytic functions}
\author[S. Sivaprasad Kumar]{S. Sivaprasad Kumar$^*$}
\address{Department of Applied Mathematics\\ Delhi Technological University\\ Delhi-110042, India}
\email{spkumar@dce.ac.in}

\author{Virendra Kumar}
\address{Department of Applied Mathematics\\ Delhi Technological University\\ Delhi-110042, India}
\email{vktmaths@yahoo.in}


\begin{document}

\begin{abstract}
In the present investigation, we derive  Fekete-Szeg\"{o} inequality
for the class $\mathcal{S}^{\alpha}_{\mathscr{L}_{g}}(\phi),$ introduced here. In addition to that, certain applications of our results are also discussed.
\end{abstract}
\footnotetext{*Corresponding author}
\maketitle

\section{Introduction}
Let $\mathcal{A}$ denote the class of functions of the form
\begin{equation}\label{eq1}
f(z)=z+\sum^{\infty}_{n=2}a_nz^n,
\end{equation}
which are analytic in the unit disc $\mathbb {U}:=\{z\in\mathbb{C}:|z|<1\}.$
Further let $\mathcal{S}$ denote the subclass of $\mathcal{A}$ consisting of univalent functions. Assume that $\phi$ is an analytic
function with positive real part in the unit disc $\mathbb {U}$ with $\phi(0)=1$ and $\phi'(0)>0$, which
maps the unit disc $\mathbb {U}$ onto a region starlike with respect to $1$ and symmetric
with respect to the real axis.

For any two analytic functions $f$ and $g$, we say that $f$ is \emph{subordinate} to $g$ or $g$ is \emph{superordinate} to $f$, denoted by $f\prec g$, if there exists a Schwarz function $w$ with $|w(z)|\leq |z|$ such that $f(z)=g(w(z)).$ If $g$ is univalent, then $f\prec g$ if and only if $f(0)=g(0)$ and $f(\mathbb{U})\subseteq g(\mathbb{U})$.  A function $p(z)=1+p_1z+p_2z^2+\ldots$ is said to be in the class $\mathcal{P}$ if $\RE{p(z)}>0.$

Let $\mathcal{S}^*(\phi)$ be the class of functions $f\in \mathcal{S}$ satisfy
$$\frac{zf'(z)}{f(z)}\prec\phi(z)\quad (z\in \mathbb {U})$$
and $\mathcal{C}(\phi)$ be the class of functions $f\in \mathcal{S}$ satisfy
$$1+\frac{zf''(z)}{f'(z)}\prec\phi(z)\quad (z\in \mathbb {U}),$$
these classes were introduced and studied by Ma and Minda \cite{minda}.
Note that $\mathcal{S}^*(\frac{1+z}{1-z})=:\mathcal{S}^*$ and $\mathcal{C}(\frac{1+z}{1-z})=:\mathcal{C}$ are the well known classes of starlike and convex functions respectively.

If $f\in \mathcal{A}$ is given by (\ref{eq1}) and $g\in \mathcal{A}$ is given by

\begin{equation}\label{g}
    g(z)=z +\sum_{n=2}^{\infty}b_n z^n,
\end{equation}
then the Hadamard product(or convolution) $f*g$ of $f$  and  $g$
 is defined by
  $$(f*g)(z):=z +\sum_{n=2}^{\infty}a_n b_n z^n=:(g*f)(z).$$

In $1933,$ M. Fekete and G. Szeg\"{o} \cite{fekete33} obtained  the sharp bounds for $|a_3-\mu a_2^2|$ as a function of real parameter $\mu$ and proved that
$$|a_2^2-\mu a_3|\leq 1+2\exp{\left(-\frac{2\mu}{1-\mu}\right)}\quad (0\leq\mu\leq1),$$
for functions belonging to the class $\mathcal{S}$. Later the problem of finding the sharp bounds for the non-linear functional $|a_3-\mu a_2^2|$ of any compact family of functions $f\in \mathcal{A}$ is known as the Fekete-Szeg\"{o} problem or inequality. In the recent years several authors have investigated the Fekete-Szeg\"{o} inequality for various subclasses of analytic functions \cite{ravi,dorina, srivastava, tuneski}.

In the present investigation, we obtain the Fekete-Szeg\"{o} inequality for functions belonging to the class $\mathcal{S}^{\alpha}_{\mathscr{L}_{g}}(\phi),$ defined here.  Applications of our main results are also discussed. In fact we generalize many earlier results in this direction \cite{koef69, minda, moorthy, dorina}.

\begin{defn} Let $f,g\in\mathcal{A}$ are respectively given by (\ref{eq1}) and (\ref{g}). We define the convolution operator $\mathscr{L}_{g}$
by $$\mathscr{L}_{g}(f(z)):=f(z)*g(z).$$
\end{defn}
We note that the operator $\mathscr{L}_g$ unifies many earlier linear operators for a suitable
choice of the function $g(z)$. Some are listed below:
\begin{enumerate}
\item If  $g(z)=\sum_{n=0}^\infty \frac{(\alpha_1)_{n}
\ldots(\alpha_l)_{n} }{ (\beta_1)_{n} \ldots (\beta_m)_{n} }
\frac{z^{n+1}}{(n)!},$ then the operator $\mathscr{L}_{g}$ coincides with  the Dziok-Srivastava~\cite{DHMS03} linear operator $\mathcal{H}^{l,m}[\alpha_1]$.
  \item  If $g(z)=\frac{z}{(1-z)^{n+1}},$ then $\mathscr{L}_g$ reduces to $D^n,$ where $D^n$ is the Ruscheweyh derivative operator \cite{rush}.
   \item For $g(z)=z+\sum_{n=2}^{\infty}\left(\frac{n+\lambda}{1+\lambda}\right)^rz^n,$
         the operator $\mathscr{L}_{g}$ reduces to the operator $I_1(r,\lambda),$
          defined by Sivaprasad {\it et al.} \cite{siva}.
  \item For $g(z)=z+\sum_{n=2}^{\infty}\frac{\Gamma{(n+1)}\Gamma{(2-\delta)}}
      {\Gamma{(n+1-\delta)}}z^n$, the operator $\mathscr{L}_{g}$ coincides with the fractional derivative operator $\Omega^\delta$, defined by Owa and Srivastava \cite{owa}.
\end{enumerate}
\begin{defn}\label{def2} Let $\alpha$ be a complex number. A function $f\in \mathcal{A}$ of the form (\ref{eq1}) is said to be in the class $\mathcal{S}^{\alpha}_{\mathscr{L}_{g}}(\phi)$ if it satisfies
\begin{equation}\label{eq2}
\Psi_{g}(f)(z)\prec\phi(z),
\end{equation}
where $$\Psi_{g}(f)(z):=1+\frac{z\mathscr{L}_{g}'(f(z))}{\mathscr{L}_{g}(f(z))}+\frac{z\mathscr{L}_{g}''(f(z))}{\mathscr{L}_{g}'(f(z))}
-\frac{(1-\alpha) z^2\mathscr{L}_{g}''(f(z))+z\mathscr{L}_{g}'(f(z))}
{(1-\alpha) z\mathscr{L}_{g}'(f(z))+\alpha \mathscr{L}_{g}(f(z))}$$ and $\mathscr{L}_{g}(f(z)):=f(z)*g(z).$
\end{defn}
\begin{rem}
    For $g(z)=\frac{z}{1-z},$ we have $\mathcal{S}^{0}_{\mathscr{L}_{g}}(\phi)=\mathcal{S}^*(\phi)$ and $\mathcal{S}^{1}_{\mathscr{L}_{g}}(\phi)=\mathcal{C}(\phi).$
  \end{rem}
\begin{rem} If we take $g(z)=z+\sum_{n=2}^{\infty}n^m z^n$, then the operator $\mathscr{L}_{g}$ reduces to the S\u{a}l\u{a}gean \cite{salagean} differential operator  $\mathcal{D}^m$  defined by \begin{eqnarray*}
   \mathcal{D}^mf(z)=  z+\sum_{n=2}^{\infty}n^m a_nz^n,  \quad m\in\{0, 1, 2, 3,\ldots\}.
\end{eqnarray*} Further, if we set $\phi(z)=\frac{1+z}{1-z}$ and  $g=z+\sum_{n=2}^{\infty}n^m z^n$ in the above Definition~\ref{def2}, then the class  $\mathcal{S}^{\alpha}_{\mathscr{L}_{g}}(\phi)$ reduces to the class $\mathcal{HS}^*_{m}(\alpha),$ recently introduced by R\u{a}ducanu \cite{dorina}.
   R\u{a}ducanu  in fact investigated the relationship property between the classes $\mathcal{HS}^*_{m}(\alpha)$ and $\mathcal{S}^*$ and obtained the Fekete-Szeg\"{o} inequality for the class $\mathcal{HS}^*_{m}(\alpha)$.
 \end{rem}
We need the following Lemmas to prove our main results:
\begin{lem}\cite{minda}\label{minda} If $p_{1}(z)=1+c_{1}z+c_{2}z^2+\ldots\in\mathcal{P}$. Then

$$|c_2-v c_1^2|\leq\left\{
  \begin{array}{ll}
    -4v+2   & \hbox{if \;  $v\leq 0$;} \\
    2    & \hbox{if \;  $0\leq v\leq 1$;} \\
    4v-2   & \hbox{if \;  $v\geq 1$.}
  \end{array}
\right.$$
When $v< 0$ or $v> 1$, equality holds if and only if $p_1(z)$ is $(1+ z)/(1 - z)$ or one of
its rotations. If  $0<v< 1$, then equality holds if and only if $p_1(z)$ is $(1+ z^2)/(1 - z^2)$
or one of its rotations. If $v=0$, equality holds if and only if
\begin{equation}\label{eq}p_1(z)=\left(\frac{1+\gamma}{2}\right)\frac{1+z}{1-z}+
\left(\frac{1-\gamma}{2}\right)\frac{1-z}{1+z} \quad (0\leq \gamma\leq1, z\in \mathbb{U})
\end{equation}
or one of its rotations. While for $v =1$, equality holds if and only if $p_1(z)$ is the reciprocal of one of the functions such that equality holds in the case of $v =0$.

 Although the above upper bound is sharp, it can be improved as follows when $0 <v< 1$:
$$|c_2-vc_1^2|+v|c_1|^2\leq 2, \quad 0 <v \leq \frac{1}{2} $$
and
$$|c_2-vc_1^2|+(1-v)|c_1|^2\leq 2,\quad \frac{1}{2} <v \leq 1 . $$
\end{lem}
\begin{lem}\cite{ravi}\label{ravi} If $p_1(z)=1+c_{1}z+c_{2}z^2+\ldots \in \mathcal{P}$. Then for any complex number $v$,
 $$|c_2-vc_1^2|\leq 2 \max\{1;|2v-1|\}$$
and the result is sharp for the functions given by
$$p_1(z) =\frac{1+z^2}{1-z^2}\;\;\text{and}\; \;p_1(z) =\frac{1+z}{1-z}.$$
\end{lem}
\begin{lem}\cite{pomer}\label{pomer} If the function $p_1(z)=1+c_{1}z+c_{2}z^2+\ldots \in \mathcal{P}$. Then
\begin{enumerate}
  \item $|c_n|\leq2$ \quad for $n\geq1$,
  \item $|c_2-\frac{1}{2}c_1^2|\leq2-\frac{|c_1|^2}{2}.$
\end{enumerate}

\end{lem}
 \section{The Fekete-Szeg\"{o} Inequality}
We begin with the following result with a coefficient estimate for the class of functions $f\in\mathcal{S}_{\mathscr{L}_{g}}^\alpha(\phi)$.
\begin{thm}\label{th}  Let  $g(z)$ be given by (\ref{g}) with $b_2,b_3$ non zero real numbers.
 Assume that $\alpha\geq0$ and $\phi(z)=1+B_1z +B_2z^2 +\cdots.$ If $f\in\mathcal{S}_{\mathscr{L}_{g}}^\alpha(\phi)$, then
\begin{equation}\label{eqn5}|a_2|\leq\frac{B_1}{(1+\alpha)|b_2|}\end{equation}
\end{thm}
and for any real number $\mu$
\begin{equation}\label{eq1.1}|a_3-\mu a_2^2|\leq\left\{
  \begin{array}{ll}
    \frac{B_1}{2(2\alpha+1)|b_3|}\left(\frac{B_2}{B_1}-
     \frac{(\alpha^2-4\alpha-1)B_1}{(1+\alpha)^2}-\frac{2\mu(2\alpha+1)B_1b_3}{(1+\alpha)^2b_2^2}\right)
        & \hbox{if \; $ \mu\leq \sigma_1$;} \\
    \frac{B_1}{2(2\alpha+1)|b_3|}
     & \hbox{if\; $ \sigma_1\leq \mu\leq \sigma_2$;} \\
    \frac{B_1}{2(2\alpha+1)|b_3|}\left(
     \frac{(\alpha^2-4\alpha-1)B_1}{(1+\alpha)^2}+\frac{2\mu(2\alpha+1)B_1b_3}{(1+\alpha)^2b_2^2}-\frac{B_2}{B_1}\right)
 & \hbox{if\; $ \mu\geq \sigma_2$,}
  \end{array}
\right. \end{equation}
where
$$\sigma_1:=\frac{(1+\alpha)^2b_2^2}{2(2\alpha+1)B_1b_3}\left(\frac{B_2}{B_1}-\frac{(\alpha^2-4\alpha-1)B_1}
{(1+\alpha)^2}-1\right) $$
and
$$\sigma_2:=\frac{(1+\alpha)^2b_2^2}{2(2\alpha+1)B_1b_3}\left(1+\frac{B_2}{B_1}-\frac{(\alpha^2-4\alpha-1)B_1}
{(1+\alpha)^2}\right).$$
The inequality (\ref{eq1.1}) is sharp.
\begin{proof} Let $f\in\mathcal{S}_{\mathscr{L}_{g}}^\alpha(\phi)$ and
\begin{eqnarray}\label{eq3}
  p(z)&=& 1+\frac{z\mathscr{L}_{g}'(f(z))}{\mathscr{L}_{g}(f(z))}+\frac{z\mathscr{L}_{g}''(f(z))}{\mathscr{L}_{g}'(f(z))}
-\frac{(1-\alpha) z^2\mathscr{L}_{g}''(f(z))+z\mathscr{L}_{g}'(f(z))}
{(1-\alpha) z\mathscr{L}_{g}'(f(z))+\alpha \mathscr{L}_{g}(f(z))} \qquad \;\\
 \nonumber  &=& 1+d_1z+d_2z^2+\cdots.
\end{eqnarray}
A simple computation shows that
$$\frac{z\mathscr{L}_{g}'(f(z))}{\mathscr{L}_{g}(f(z))}=1+a_2b_2z+[2a_3b_3-a_2^2b_2^2]z^2+\ldots, $$
$$1+\frac{z\mathscr{L}_{g}''(f(z))}{\mathscr{L}_{g}'(f(z))}=1+2a_2b_2z+[6a_3b_3-4a_2^2b_2^2]z^2+\cdots $$
and
$$\frac{(1-\alpha) z^2\mathscr{L}_{g}''(f(z))+z\mathscr{L}_{g}'(f(z))}
{(1-\alpha) z\mathscr{L}_{g}'(f(z))+\alpha \mathscr{L}_{g}(f(z))}=1+(2-\alpha)a_2b_2z+[(6-4\alpha)a_3b_3-(\alpha-2)^2a_2^2b_2^2]z^2+\cdots. $$
Substituting these values in (\ref{eq3}), we have
\begin{equation}\label{eq4}d_1=(1+\alpha) a_2b_2\end{equation}
 and \begin{equation}\label{eq5}d_2=2(2\alpha+1)a_3b_3+(\alpha^2-4\alpha-1)a_2^2b_2^2.\end{equation}
Since $\phi$ is univalent and $p\prec\phi$, the function $p_1(z)$ defined by
\begin{equation}\label{eqn}p_1(z)=\frac{1+\phi^{-1}(p(z))}{1-\phi^{-1}(p(z))}=1+ c_1z + c_2z^2
+ \ldots, \end{equation}
 is analytic with positive real part in the unit disc $\mathbb{U}.$ Further from (\ref{eqn}),
we have
\begin{eqnarray*}
  p(z) &=& \phi\left(\frac{p_1(z)-1}{p_1(z)+1}\right) \\
  &=& \phi\left(\frac{c_1z + c_2z^2+\ldots}{2+c_1z + c_2z^2+\ldots}\right) \\
   &=& 1+\frac{1}{2}B_1c_1z+\left[\frac{1}{2}B_1(c_2-\frac{1}{2}c_1^2)+\frac{1}{4}B_2c_1^2\right]z^2
   +\cdots.
\end{eqnarray*}
Thus, we have
\begin{equation}\label{eqn1}d_1=\frac{1}{2}B_1c_1\end{equation} and \begin{equation}\label{eqn2}d_2=\frac{1}{2}B_1\left(c_2-\frac{1}{2}c_1^2\right)+\frac{1}{4}B_2c_1^2.\end{equation}
From (\ref{eq4}) and (\ref{eqn1}), we have
\begin{equation}\label{eqn3}a_2=\frac{B_1c_1}{2(1+\alpha)b_2}. \end{equation}
Similarly from (\ref{eq5}) and (\ref{eqn2}), we obtain
\begin{equation}\label{eqn4}a_3=\frac{[2B_1(c_2-\frac{1}{2}c_1^2)+B_2c_1^2](1+\alpha)^2-(\alpha^2-4\alpha-1)B_2^2c_1^2}
{8(2\alpha+1)(1+\alpha)^2b_3}.\end{equation}
The inequality (\ref{eqn5}) now follows from (\ref{eqn3}) and the first part of Lemma~\ref{pomer}.

By using (\ref{eqn3}) and (\ref{eqn4}), we have

 \begin{equation}\label{a1}
    a_3-\mu a_2^2 = \frac{B_1}{4(2\alpha+1)b_3}[c_2-vc_1^2],
\end{equation}
where $$v:=\frac{1}{2}\left[1-\frac{B_2}{B_1}+\frac{(\alpha^2-4\alpha-1)B_1}{(1+\alpha)^2}+
\frac{2\mu(2\alpha+1)B_1b_3}{(1+\alpha)^2b_2^2}\right].$$
If $\mu\leq\sigma_1$, then an application of Lemma~\ref{minda} gives
$$|a_3-\mu a_2^2|\leq\frac{B_1}{2(2\alpha+1)|b_3|}\left(\frac{B_2}{B_1}-
     \frac{(\alpha^2-4\alpha-1)B_1}{(1+\alpha)^2}-\frac{\mu(2\alpha+1)B_1b_3}
     {(1+\alpha)^2b_2^2}\right),$$
which is the first part of assertion (\ref{eq1.1}).

Next, if $\mu\geq\sigma_2$, then by applying Lemma~\ref{minda}, we can write
$$|a_3-\mu a_2^2|\leq\frac{B_1}{2(2\alpha+1)|b_3|}\left(
     \frac{(\alpha^2-4\alpha-1)B_1}{(1+\alpha)^2}+\frac{\mu(2\alpha+1)B_1b_3}
{(1+\alpha)^2b_2^2}-\frac{B_2}{B_1}\right),$$
which is the third part of assertion (\ref{eq1.1}).

If $\mu=\sigma_1$, then equality holds if and only if $p_1(z)$ is given by (\ref{eq}) or
one of its rotations.

If $\mu=\sigma_2$, then
$$\frac{1}{2}\left[1-\frac{B_2}{B_1}+\frac{(\alpha^2-4\alpha-1)B_1}{(1+\alpha)^2}+
\frac{2\mu(2\alpha+1)B_1b_3}
{(1+\alpha)^2b_2^2}\right]=1.$$
 Therefore,
 $$\frac{1}{p_1(z)}=\left(\frac{1+\gamma}{2}\right)\frac{1+z}{1-z}+
\left(\frac{1-\gamma}{2}\right)\frac{1-z}{1+z} \quad (0<\gamma<1, z\in \mathbb{U}).$$
Finally, we see that
 $$ a_3-\mu a_2^2 = \frac{B_1}{4(2\alpha+1)b_3}\left[c_2-\frac{c_1^2}{2}\left(1-\frac{B_2}{B_1}+\frac{(\alpha^2-4\alpha-1)B_1}{(1+\alpha)^2}+
\frac{2\mu(2\alpha+1)B_1b_3}{(1+\alpha)^2b_2^2}\right)\right].$$
Therefore using Lemma~\ref{minda}, we get
$$|a_3-\mu a_2^2|\leq \frac{B_1}{2(2\alpha+1)|b_3|}\quad
(\sigma_1\leq\mu\leq\sigma_2).$$
If $\sigma_1<\mu<\sigma_2$, then we have
$$p_1(z)=\frac{1+\lambda z^2}{1-\lambda z^2} \quad (0\leq\lambda \leq1).$$
By an application of Lemma~\ref{minda}, we obtain our result. To show that the inequality (\ref{eq1.1}) is sharp,
we define the functions
 $K^{\phi_n}\; (n=2,3,4,...)$ by
$$\Psi_{g}(K^{\phi_n})(z)=\phi(z^{n-1})\quad (K^{\phi_n}(0)=0=(K^{\phi_n})'(0)-1)$$
and the functions $G^\gamma$ and $H^\gamma \; (0\leq\gamma\leq1)$ by
$$ \Psi_{g}(G^\gamma)(z)=\phi\left(\frac{z(z+\gamma)}{1+\gamma z}\right)\quad(
G^\gamma(0)=0=(G^\gamma)'(0)-1)$$
and
$$\Psi_{g}(H^\gamma)(z)=\phi\left(-\frac{z(z+\gamma)}{1+\gamma z}\right)\quad(
H^\gamma(0)=0=(H^\gamma)'(0)-1).$$
It is clear that the functions $K^{\phi_n}\; (n=2,3,4,...), G^\gamma $ and $H^\gamma \;
(0\leq\gamma\leq1)$ are in the class $S^{\alpha}_{\mathscr{L}_{g}}(\phi).$
In either cases $\mu<\sigma_1$ or $\mu>\sigma_2$, the equality holds if and only if $f$ is $K^{\phi_2}$ or one of its rotations. When $\sigma_1<\mu<\sigma_2$ the equality occurs
if and only if $f$ is $K^{\phi_3}$ or one of its rotations. If $\mu=\sigma_1$, then the equality holds if and only if $f$ is $G^\gamma$ or one of its rotations. If $\mu=\sigma_2$, then the equality holds if and only if $f$ is $H^\gamma$ or one of its rotations.
\end{proof}

\begin{rem}\label{cor1} Using Lemma~\ref{minda} the result can be improved when $\sigma_1\leq\mu\leq\sigma_2$ as follows:

Let
$$\sigma_3:=\frac{(1+\alpha)^2b_2^2}{2(2\alpha+1)B_1b_3}\left(\frac{B_2}{B_1}-
\frac{(\alpha^2-4\alpha-1)B_1}{(1+\alpha)^2}\right).$$
If $\sigma_1\leq\mu\leq\sigma_3$, then
\begin{align*}
   | a_3-\mu a_2^2|+\frac{(1+\alpha)^2b_2^2}{2(2\alpha+1)B_1|b_3|}&\left(1-\frac{B_2}{B_1}
+\frac{(\alpha^2-4\alpha-1)B_1}{(1+\alpha)^2}+\frac{2\mu(2\alpha+1)B_1b_3}{(1+\alpha)^2b_2^2}\right)
|a_2|^2\\
&\leq\frac{B_1}{2(2\alpha+1)|b_3|}
\end{align*}
and if $\sigma_3\leq\mu\leq\sigma_2$, then
\begin{align*}
   | a_3-\mu a_2^2|+\frac{(1+\alpha)^2b_2^2}{2(2\alpha+1)B_1|b_3|}&\left(1+\frac{B_2}{B_1}
-\frac{(\alpha^2-4\alpha-1)B_1}{(1+\alpha)^2}-\frac{2\mu(2\alpha+1)B_1b_3}
{(1+\alpha)^2b_2^2}\right)|a_2|^2\\ & \leq\frac{B_1}{2(2\alpha+1)|b_3|}.
\end{align*}

\begin{proof}
For the values  of $\sigma_1\leq\mu\leq\sigma_3$, we have
\begin{align*}
| a_3-\mu a_2^2|+\frac{(1+\alpha)^2b_2^2}{2(2\alpha+1)B_1|b_3|}&\left(1-\frac{B_2}{B_1}
+\frac{(\alpha^2-4\alpha-1)B_1}{(1+\alpha)^2}+\frac{2\mu(2\alpha+1)B_1b_3}{(1+\alpha)^2b_2^2}\right)
|a_2|^2\\&=| a_3-\mu a_2^2|+(\mu-\sigma_1)|a_2|^2
\\&=\frac{B_1}{4(2\alpha+1)|b_3|}\left[|c_2-vc_1^2|+v|c_1|^2\right]\\
&\leq \frac{B_1}{2(2\alpha+1)|b_3|}.
\end{align*}
Similarly, if $\sigma_3\leq\mu\leq\sigma_2$, then

\begin{align*}
  | a_3-\mu a_2^2|+\frac{(1+\alpha)^2b_2^2}{2(2\alpha+1)B_1|b_3|}&\left(1+\frac{B_2}{B_1}
-\frac{(\alpha^2-4\alpha-1)B_1}{(1+\alpha)^2}-\frac{2\mu(2\alpha+1)B_1b_3}
{(1+\alpha)^2b_2^2}\right)|a_2|^2\\ &= | a_3-\mu a_2^2|+(\sigma_2-\mu)|a_2|^2 \\
   &= \frac{B_1}{4(2\alpha+1)|b_3|}\left[|c_2-vc_1^2|+(1-v)|c_1|^2\right]\\
  &\leq \frac{B_1}{2(2\alpha+1)|b_3|}.
\end{align*}
Thus the proof is complete.
\end{proof}\end{rem}
\begin{rem}
If we set $\alpha=1$ and $g(z)=z/(1-z)$ in Theorem~\ref{th}, then
we have the result \cite[Theorem 3]{minda} of Ma and Minda.
\end{rem}
\begin{rem} By setting $\alpha=0$ and $g(z)=z/(1-z)$ in Theorem~\ref{th}, we obtain the
result of Murugusundaramoorthy {\it et al.} \cite[Corollary 2.2]{moorthy}.
\end{rem}
Using Lemma~\ref{ravi} and equation (\ref{a1}), we deduce the following:
\begin{thm}\label{th2} Let  $g(z)$ be given by (\ref{g}) with $b_2,b_3$ non zero real numbers.
 Assume that $\alpha\geq0$ and $\phi(z)=1+B_1z +B_2z^2 +\cdots.$
If $f\in\mathcal{S}_{\mathscr{L}_{g}}^\alpha(\phi)$, then for any
complex number $\mu$
$$|a_3-\mu a_2^2|\leq \frac{B_1}{2(2\alpha+1)|b_3|} \max \set{1; \left|\frac{2\mu (2\alpha+1)B_1b_3}{(1+\alpha)^2b_2^2}+\frac{(\alpha^2-4\alpha-1)B_1}{(1+\alpha)^2}-\frac{B_2}{B_1}\right|} .$$
\end{thm}
\section{Applications}
A few applications of our main results are discussed here.
\begin{defn}\cite{srivastava1}\label{df1}
Let $f(z)$ be an analytic function in a simply connected region
of the complex plane containing the origin. The fractional derivative of order $\delta$ is defined by
$$D^\delta_zf(z)=\frac{1}{\Gamma(1-\delta)}\frac{d}{dz}\int^{z}_0\frac{f(t)}{(z-t)^\delta}dt \quad (0\leq\delta<1),$$
where the multiplicity of $(z-t)^\delta$ is removed by requiring that $\log(z-t)$ is real
for $(z-t)>0$.
\end{defn}
Using the above Definition~\ref{df1} and its extensions involving fractional derivatives and
fractional integrals, Owa and Srivastava \cite{owa} introduced the operator $\Omega^\delta:
\mathcal{A}\rightarrow \mathcal{A}$ defined by
$$(\Omega^\delta f)(z)=\Gamma{(2-\delta)}z^\delta D^\delta_zf(z), \; \delta\neq 2,3,4\cdots.$$
If we take $g(z)=z+\sum^{\infty}_{n=2}\frac{\Gamma{(n+1)}
\Gamma{(2-\delta)}}{\Gamma{(n+1-\delta)}}z^n$ in Theorem \ref{th}, then we obtain the
following:
\begin{cor} Let $\alpha\geq0$ and $\phi(z)=1+B_1z +B_2z^2 +\ldots\;.$ If $f\in\mathcal{S}_{g}^\alpha(\phi)$, then
$$|a_2|\leq \frac{(2-\delta)B_1}{2(1+\alpha)}$$
and for any real number $\mu$
$$|a_3-\mu a_2^2|\leq
\left\{
  \begin{array}{ll}
    \frac{(2-\delta)(3-\delta)B_1}{12(2\alpha+1)}\left(\frac{B_2}{B_1}-
     \frac{(\alpha^2-4\alpha-1)B_1}{(1+\alpha)^2}-\frac{3\mu(2\alpha+1)(2-\delta)B_1}
      {(1+\alpha)^2(3-\delta)}\right)
       & \hbox{if\; $\mu\leq \sigma_1$;} \\
    \frac{(2-\delta)(3-\delta)B_1}{12(2\alpha+1)}
    & \hbox{if\;$\sigma_1\leq \mu\leq \sigma_2$;} \\
    \frac{(2-\delta)(3-\delta)B_1}{12(2\alpha+1)}\left(
     \frac{(\alpha^2-4\alpha-1)B_1}{(1+\alpha)^2}+\frac{3\mu(2\alpha+1)(2-\delta)B_1}
    {(1+\alpha)^2(3-\delta)}-\frac{B_2}{B_1}\right)  & \hbox{if\;$\mu\geq \sigma_2$,}
  \end{array}
\right.$$
where
$$\sigma_1:=\frac{(1+\alpha)^2(3-\delta)}{3(2-\delta)(2\alpha+1)B_1}\left(\frac{B_2}{B_1}
-\frac{(\alpha^2-4\alpha-1)B_1}{(1+\alpha)^2}-1\right) $$
and
$$\sigma_2:=\frac{(1+\alpha)^2(3-\delta)}{(2-\delta)(2\alpha+1)B_1}\left(1+\frac{B_2}{B_1}-
\frac{(\alpha^2-4\alpha-1)B_1}{(1+\alpha)^2}\right).$$
The result is sharp.
\end{cor}
From Theorem~\ref{th} and Remark ~\ref{cor1}, we deduce the following:
\begin{cor} \label{cor2}  Let  $g(z)$ be given by (\ref{g}) with $b_2,b_3$ non zero real numbers. Assume that $\alpha\geq0$ and $-1\leq D<C\leq1$.
 If $f\in\mathcal{S}_{\mathscr{L}_{g}}^\alpha(\frac{1+Cz}{1+Dz})$, then
$$|a_2|\leq\frac{C-D}{(1+\alpha)|b_2|}$$
 and for any real number $\mu$
$$|a_3-\mu a_2^2|\leq\left\{
  \begin{array}{ll}
    \frac{D-C}{2(2\alpha+1)|b_3|}\left(D+
     \frac{(\alpha^2-4\alpha-1)(C-D)}{(1+\alpha)^2}+\frac{2\mu(2\alpha+1)(C-D)b_3}{(1+\alpha)^2b_2^2}\right)
        & \hbox{if \; $\mu\leq \sigma_1$;} \\
    \frac{C-D}{2(2\alpha+1)|b_3|}
     & \hbox{if \;$\sigma_1\leq \mu\leq \sigma_2$;} \\
    \frac{C-D}{2(2\alpha+1)|b_3|}\left(D+
     \frac{(\alpha^2-4\alpha-1)(C-D)}{(1+\alpha)^2}+\frac{2\mu(2\alpha+1)(C-D)b_3}{(1+\alpha)^2b_2^2}\right)
 & \hbox{if \;$\mu\geq \sigma_2$,}
  \end{array}
\right. $$
where
$$\sigma_1:=\frac{(1+\alpha)^2b_2^2}{2(2\alpha+1)(D-C)b_3}\left(1+D+\frac{(\alpha^2-4\alpha-1)(C-D)}
{(1+\alpha)^2}\right) $$
and
$$\sigma_2:=\frac{(1+\alpha)^2b_2^2}{2(2\alpha+1)(C-D)b_3}\left(1-D-\frac{(\alpha^2-4\alpha-1)(C-D)}
{(1+\alpha)^2}\right).$$ The result is sharp.\end{cor}
\begin{rem}
  The result can be improved when $\sigma_1\leq\mu\leq\sigma_2$ as follows:

If $\sigma_1\leq\mu\leq\sigma_3$, then
{\small \begin{align*}
  | a_3-\mu a_2^2|+\frac{(1+\alpha)^2b_2^2}{2(2\alpha+1)(C-D)|b_3|}&\left(1+D
+\frac{(\alpha^2-4\alpha-1)(C-D)}{(1+\alpha)^2}+
\frac{2\mu(2\alpha+1)(C-D)b_3}{(1+\alpha)^2b_2^2}\right)|a_2|^2\\&
\leq\frac{C-D}{2(2\alpha+1)|b_3|}
\end{align*}}
and if $\sigma_3\leq\mu\leq\sigma_2$, then
{\small \begin{align*}
   | a_3-\mu a_2^2|+\frac{(1+\alpha)^2b_2^2}{2(2\alpha+1)(C-D)|b_3|}&\left(1-D
-\frac{(\alpha^2-4\alpha-1)(C-D)}{(1+\alpha)^2}
-\frac{2\mu(2\alpha+1)(C-D)b_3}{(1+\alpha)^2b_2^2}\right)|a_2|^2\\&
 \leq\frac{C-D}{2(2\alpha+1)|b_3|},
\end{align*}}
where
$$\sigma_3:=\frac{(1+\alpha)^2b_2^2}{2(2\alpha+1)(D-C)b_3}\left(D+\frac{(\alpha^2-4\alpha-1)(C-D)}
{(1+\alpha)^2}\right).$$
\end{rem}

By taking $D=-1$ and $C=1$ in the above Corollary~\ref{cor2}, we obtain the following:
\begin{exm}\label{ex3} Let $\alpha\geq0$ and $g(z)$ be given by (\ref{g}) with $b_2,b_3$ non zero real numbers.
 If $f\in\mathcal{S}_{\mathscr{L}_{g}}^\alpha(\frac{1+z}{1-z})$, then
$$|a_2|\leq\frac{2}{(1+\alpha)|b_2|}$$
 and for any real number $\mu$
$$|a_3-\mu a_2^2|\leq\left\{
  \begin{array}{ll}
    \frac{1}{(1+\alpha)^2|b_3|}\left(
     \frac{3+10\alpha-\alpha^2}{2\alpha+1}-\frac{4\mu b_3}{b_2^2}\right)
        & \hbox{if\; $\mu\leq \sigma_1$;} \\
    \frac{1}{(2\alpha+1)|b_3|}
     & \hbox{if\; $\sigma_1\leq \mu\leq \sigma_2$;} \\
    \frac{1}{(1+\alpha)^2|b_3|}\left(
     \frac{\alpha^2-10\alpha-3}{2\alpha+1}+\frac{4\mu b_3}{b_2^2}\right)
 & \hbox{if\; $\mu\geq \sigma_2$,}
  \end{array}
\right. $$
where
$$\sigma_1:=\frac{(1+4\alpha-\alpha^2)b_2^2}{2(2\alpha+1)b_3}\;\; {\text {and}}\;\;
\sigma_2:=\frac{(3\alpha+1)b_2^2}{(2\alpha+1)b_3}.$$
The result can be improved when $\sigma_1\leq\mu\leq\sigma_2$ as follows:

If $\sigma_1\leq\mu\leq\sigma_3$, then
$$ | a_3-\mu a_2^2|+\frac{b_2^2}{2|b_3|}\left(
\frac{\alpha^2-4\alpha-1}{2\alpha+1}+\frac{2\mu b_3}{b_2^2}\right)|a_2|^2
\leq\frac{1}{(2\alpha+1)|b_3|}
$$
and if $\sigma_3\leq\mu\leq\sigma_2$, then
$$ | a_3-\mu a_2^2|+\frac{b_2^2}{|b_3|}\left(\frac{3\alpha+1}{2\alpha+1}-\frac{\mu b_3}{b_2^2}\right)|a_2|^2 \leq \frac{1}{(2\alpha+1)|b_3|}, $$
where $$\sigma_3:=\frac{(3+10\alpha-\alpha^2)b_2^2}{4(2\alpha+1)b_3}.$$
The result is sharp.
\end{exm}
\begin{rem}We obtain the result of R\u{a}ducanu \cite[Theorem 2]{dorina} by taking
$$
   g(z)=z+\sum_{n=2}^{\infty}n^mz^n \quad(m\in\{0, 1, 2, 3,\ldots\})
 $$
in the above Example~\ref{ex3}.
\end{rem}
\begin{rem}
Setting $\alpha=0$ and $g(z)=\frac{z}{1-z}$ in the Example~\ref{ex3}, we obtain the following result \cite{srivastava}:

 If $f\in\mathcal{S}^*$, then
$$|a_2|\leq 2$$
and for any real number $\mu$
$$|a_3-\mu a_2^2|\leq\left\{
  \begin{array}{ll}
    3-4\mu
        & \hbox{ if $\;\mu\leq \frac{1}{2}$;} \\
    1
     & \hbox{ if $\;\frac{1}{2}\leq \mu\leq 1$;} \\
    4\mu-3
 & \hbox{ if $\;\mu\geq 1$.}
  \end{array}
\right. $$
The result can be improved when $\frac{1}{2}\leq\mu\leq 1$ as follows:

If $\frac{1}{2}\leq\mu\leq \frac{3}{4}$, then
$$ | a_3-\mu a_2^2|+\frac{1}{2}(2\mu-1)|a_2|^2\leq 1$$
and
if $\frac{3}{4}\leq\mu\leq 1$, then
$$ | a_3-\mu a_2^2|+(1-\mu)|a_2|^2 \leq 1.$$

Setting $\alpha=1$ and $g(z)=\frac{z}{1-z}$ in the Example~\ref{ex3}, we have the following result:

Let $f\in\mathcal{C}$, then
$$|a_2|\leq 1$$
 and for any real number $\mu$
$$|a_3-\mu a_2^2|\leq\left\{
  \begin{array}{ll}
    1-\mu
        & \hbox{ if $\;\mu\leq \frac{2}{3}$;} \\
    \frac{1}{3}
     & \hbox{ if $\;\frac{2}{3}\leq \mu\leq \frac{4}{3}$;} \\
    \mu-1
 & \hbox{ if $\;\mu\geq\frac{4}{3}$.}
  \end{array}
\right. $$
The result can be improved when $\frac{2}{3}\leq\mu\leq\frac{4}{3}$ as follows:

If $\frac{2}{3}\leq\mu\leq 1$, then
$$ | a_3-\mu a_2^2|+\frac{1}{3}(3\mu-2)|a_2|^2
\leq\frac{1}{3}
$$
and if $1\leq\mu\leq\frac{4}{3}$, then
$$ | a_3-\mu a_2^2|+\frac{1}{3}(4-3\mu)|a_2|^2 \leq \frac{1}{3}.$$
\end{rem}
Taking $\phi(z)=\frac{1+Cz}{1+Dz}, -1\leq D<C\leq1$, in Theorem~\ref{th2}, we deduce the following:
\begin{cor}\label{cor3} Let $\alpha\geq0$ and $g(z)$ be given by (\ref{g}) with $b_2,b_3$ non zero real numbers. If $f\in\mathcal{S}_{\mathscr{L}_{g}}^\alpha(\frac{1+Cz}{1+Dz})$, then for any complex number $\mu$
$$|a_3-\mu a_2^2|\leq \frac{C-D}{2(2\alpha+1)|b_3|} \max \set{1; \left|\frac{2\mu (2\alpha+1)(C-D)b_3}{(1+\alpha)^2b_2^2}+\frac{(\alpha^2-4\alpha-1)(C-D)}{(1+\alpha)^2}+D\right|}. $$
\end{cor}
\begin{rem}\label{rem1}
If we take $g(z)=z+\sum_{n=2}^{\infty}n^mz^n$, $D=-1 \;\text{and}\; C=1$ in the above Corollary~\ref{cor3}, we have the following result \cite[Theorem 3]{dorina} of R\u{a}ducanu:

 Let $\alpha\geq0$. If $f\in\mathcal{HS}^*_{m}(\alpha)$, then for any complex number $\mu$
$$|a_3-\mu a_2^2|\leq \frac{1}{3^m(1+2\alpha)} \max \set{1; \frac{\left|2^{2m-1}(\alpha^2-10\alpha-3)
+2.3^m(1+2\alpha)\mu\right|}{2^{2m-1}(1+\alpha)^2}}.$$
\end{rem}
 \begin{rem}If we set $D=-1, C=1$ and $g(z)=\frac{z}{1-z}$ in Corollary \ref{cor3}, then for $\alpha=0$, we have the following result \cite[Theorem 1]{koef69}(see also \cite{srivastava}):

 Let $f\in \mathcal{S}^*$.
 Then for any complex number $\mu$
$$|a_3-\mu a_2^2|\leq \max \set{1; |4\mu-3|}.$$

 Setting $\alpha=1, D=-1, C=1$ and $g(z)=\frac{z}{1-z}$ in Corollary \ref{cor3}, we obtain the following result \cite[Corollary 1]{koef69} due to Keogh and Merkes:

Let $f\in\mathcal{C}$, then for any complex number $\mu$
$$|a_3-\mu a_2^2|\leq  \max \set{\frac{1}{3}; \left|\mu-1\right|}.$$

 \end{rem}

\end{document}